\documentclass[11pt]{article}

\usepackage{amsfonts, amsmath, amssymb, amsthm}
\usepackage[colorlinks=true, citecolor=blue]{hyperref}
\usepackage[margin=2.54cm]{geometry}
\usepackage{mathtools}
\usepackage{enumerate}
\usepackage{verbatim}
\usepackage{tikz}
\usepackage{microtype}
\usepackage{color, tikz, mathtools}
\usepackage{fontenc}
\usepackage[utf8]{inputenc}

\usepackage{tcolorbox}
\usepackage{mdframed}

\usetikzlibrary{positioning}
\usetikzlibrary{shapes}
\usetikzlibrary{decorations.pathreplacing}
\usepackage{enumerate, enumitem}

\newcommand{\RR}{\mathbb{R}}

\renewcommand{\epsilon}{\varepsilon}

\newtheorem{theorem}{Theorem}

\theoremstyle{definition}

\newcommand{\remove}[1]{}
\begin{document}

\title{
  No Infinite $(p,q)$-Theorem for Piercing Compact Convex Sets in $\mathbb{R}^3$ with Lines
}

\author{
	Sutanoya Chakraborty \footnote{Indian Statistical Institute, Kolkata, India}
	\and 
	Arijit Ghosh \footnotemark[1]
}

\date{}

\maketitle

\begin{abstract}
An infinite $(p,q)$-theorem, or an $(\aleph_0,q)$-theorem, involving two families $\mathcal{F}$ and $\mathcal{G}$ of sets, states that if in every infinite subset of $\mathcal{F}$, there are $q$ sets that are intersected by some set in $\mathcal{G}$, then there is a finite set $S_{\mathcal{F}}\subseteq\mathcal{G}$ such that for every $C\in\mathcal{F}$, there is a $B\in S_{\mathcal{F}}$ with $C\cap B\neq\emptyset$.
We provide an example demonstrating that there is no $(\aleph_0,q)$-theorem for piercing compact convex sets in $\RR^3$ with lines by constructing a family $\mathcal{F}$ of compact convex sets such that it does not have a finite line transversal, but for any $t\in\mathbb{N}$, every infinite subset of $\mathcal{F}$ contains $t$ sets that are pierced by a line.
\end{abstract}

\section{Introduction}
Given a collection $\mathcal{F}$ of sets in $\mathbb{R}^d$ and $k\in\{0,\dots,d-1\}$, a {\em $k$-transversal} of $\mathcal{F}$ is 
a family $\mathcal{K}$ consisting of $k$-flats such that for every $C\in\mathcal{F}$, there is a $K_C\in\mathcal{K}$ with $C\cap K_C\neq\emptyset$.
We say that $\mathcal{K}$ is a point transversal if $k=0$, and a line transversal if $k=1$.
Given $A,B\subseteq\mathbb{R}^d$, we say that $A$ {\em pierces} $B$ if $A\cap B\neq\emptyset$.

The $(p,q)$-theorem for piercing compact convex sets with points is a fundamental result in combinatorial geometry, proved by Alon and Kleitman \cite{alon1992piercing}.

\begin{theorem}[Alon-Kleitman $(p,q)$-theorem]
  Let $d\in\mathbb{N}$, $p,q\in\mathbb{N}$ with $p\geq q\geq d+1$, and $\mathcal{F}$ be a family of compact convex sets in $\mathbb{R}^d$ such that every 
  subset of $\mathcal{F}$ of size $p$ contains $q$ sets whose intersection is nonempty.
  Then there is an integer $c=c(p,q,d)$ such that $\mathcal{F}$ has a point transversal of size at most $c$.
\end{theorem}

Alon and Kalai \cite{alon1995bounding} proved a $(p,q)$-theorem for piercing compact connected sets with hyperplanes.

Alon, Kalai, Matou\v{s}ek and Meshulam \cite{alon2002transversal} proved that in an abstract setting, the appropriate fractional Helly properties imply a $(p,q)$-theorem. 
In the same paper, they demonstrated, by proving the following theorem, that there cannot be any {\em general} $(p,q)$-theorem for piercing compact convex sets in $\mathbb{R}^d$ by $k$-flats, without additional assumptions on the sets or the value of $k$. 
\begin{theorem}[Impossibility of a $(p,q)$-theorem for piercing convex compact sets in $\mathbb{R}^3$ with lines]
  For any $n,m,t\in\mathbb{N}$, there is a finite family $\mathcal{F}$ of compact convex sets in $\mathbb{R}^3$ with size greater than $n$,
  such that there is a line piercing every $t$ sets of $\mathcal{F}$, but any line transversal of $\mathcal{F}$ has 
  size at least $m$.
\end{theorem}

Holmsen and Matou\v{s}ek \cite{holmsen2004no} showed that we cannot have a Helly theorem for piercing translates of compact convex sets in $\mathbb{R}^3$ with lines
with the following theorem:

\begin{theorem}[Impossibility of a Helly theorem for piercing translates of convex compact sets in $\mathbb{R}^3$ with lines]
  For any $n\in\mathbb{N}$, $n>2$, there is a compact convex set $K$ and a family $\mathcal{F}$ consisting of 
  $n$ translates of $K$ in $\mathbb{R}^3$,
  such that there is a line piercing every $n-1$ sets of $\mathcal{F}$, but no single line pierces every set of $\mathcal{F}$.
\end{theorem}

Using a lifting technique by Bukh \cite{cheong2024some}, these results can be extended to $\mathbb{R}^d$.

Keller and Perles \cite{keller2022} introduced the notion of the $(\aleph_0,q)$-property.
A family $\mathcal{F}$ of sets in $\mathbb{R}^d$ satisfies the $(\aleph_0,q)$-property with respect to 
$k$-flats, $k\in\{0,\dots,d-1\}$, if every infinite subset of $\mathcal{F}$ contains $q$ sets that are 
pierced by a $k$-flat.
Note that if $\mathcal{F}$ satisfies the $(p,q)$-property, it also satisfies the $(\aleph_0,q)$-property,
but the converse is not necessarily true.

Keller and Perles \cite{keller2023aleph_0} proved the following theorem:
\begin{theorem}[$(\aleph_0,k+2)$-Theorem for piercing compact near-balls with $k$-flats]
  Let $d\in\mathbb{N}$, $k\in\{0,\dots,d-1\}$, and $\mathcal{F}$ be a family of compact sets in $\mathbb{R}^d$ that 
  are {\em near-balls} with parameter $\rho\geq 1$, i.e., every $B\in\mathcal{F}$ satisfies $\rho r_B\geq R_B$,
  $\rho + r_B \geq R_B$, where $0<r_B\geq R_B$ and there is a $b\in B$ such that $\overline{B}(b,r_B)\subseteq B\subseteq 
  \overline{B}(b,R_B)$.
  Then, if $\mathcal{F}$ satisfies the $(\aleph_0,k+2)$-property with respect to $k$-flats, then there is a finite collection 
  of $k$-flats that is a transversal of $\mathcal{F}$.
\end{theorem}

Jung and P{\'a}lv{\"o}lgyi \cite{jung2024note} proved that the existence of a fractional Helly theorem and a $(p,q)$-theorem implies the corresponding
$(\aleph_0,q)$-theorem.
In particular, they proved the following theorem:
\begin{theorem}[$(\aleph_0,d+1)$-Theorem for piercing compact convex sets with points and hyperplanes]
  Let $d\in\mathbb{N}$ and $\mathcal{F}$ be a family of compact convex sets in $\mathbb{R}^d$ that satisfies 
  the $(\aleph_0,d+1)$-property with respect to points (hyperplanes).
  Then $\mathcal{F}$ has a finite-sized point (hyperplane) transversal.
\end{theorem}

Chakraborty, Ghosh and Nandi \cite{chakraborty2025finitektransversalsinfinitefamilies} proved the following theorem: 
\begin{theorem}[$(\aleph_0,k+2)$-Theorem for piercing compact convex $\rho$-fat sets with $k$-flats]
  Let $d\in\mathbb{N}$, $k\in\{0,\dots,d-1\}$, and $\mathcal{F}$ be a family of compact convex sets in $\mathbb{R}^d$ that 
  are {\em $\rho$-fat} for some $\rho\geq 1$, i.e., every $B\in\mathcal{F}$ satisfies $\rho r_B\geq R_B$,
  where $0<r_B\leq R_B$ and there is a $b\in B$ such that $\overline{B}(b,r_B)\subseteq B\subseteq 
  \overline{B}(b,R_B)$.
  Then, if $\mathcal{F}$ satisfies the $(\aleph_0,k+2)$-property with respect to $k$-flats, then there is a finite collection 
  of $k$-flats that is a transversal of $\mathcal{F}$.
\end{theorem}

Chakraborty, Ghosh and Nandi \cite{chakraborty2025stabbing} also proved an $(\aleph_0,k+2)$-theorem for piercing axis-parallel boxes with 
axis-parallel $k$-flats, of which Jung and P{\'a}lv{\"o}lgyi \cite{jung2024note} later gave an alternative proof.

It is worth noting that in an $(\aleph_0,q)$-theorem, no bounds on the transversal size need to be given; one only needs to show that there is a finite transversal.
This is unlike $(p,q)$-theorems, where there is a universal bound on the transversal size, depending on $p,q$ and the dimension size.
In fact, one can easily see that such a bound cannot exist for any of the classes of sets mentioned in the $(\aleph_0,q)$-theorems stated above.

\section{Our result}
Chakraborty, Ghosh and Nandi \cite{chakraborty2025finitektransversalsinfinitefamilies} showed, with an example, that the existence of an $(\aleph_0,q)$-theorem does not
imply the existence of a corresponding $(p,q)$-theorem. 
Therefore, it is possible to have an $(\aleph_0,q)$-theorem even when the corresponding $(p,q)$-theorem does not exist.

We demonstrate, with the following example, that it is not possible to have an $(\aleph_0,q)$-theorem for piercing 
compact convex sets in $\mathbb{R}^d$ with $k$-flats in general, without additional assumptions on the sets or the value of $k$.

\begin{theorem}\label{exmpl:lines_R3}
    There is a collection $\mathcal{F}$ of compact convex sets in $\mathbb{R}^3$, such that for any $t\in\mathbb{N}$, 
    every infinite subset of $\mathcal{F}$ contains $t$ sets that are pierced by a line, but $\mathcal{F}$ does not have a finite line transversal.
\end{theorem}

Like Alon et al \cite{alon2002transversal} and Holmsen and Matou\v{s}ek \cite{holmsen2004no}, this construction too depends on the geometry of the hyperbolic paraboloid.
Unlike the families of sets constructed in \cite{alon2002transversal} and \cite{holmsen2004no} which are finite collections that depend on the value of $t$, the family $\mathcal{F}$ that we construct here works for 
any value of $t\in\mathbb{N}$, i.e., $\mathcal{F}$ does not have a finite line transversal, but for 
any $t\in\mathbb{N}$, in every infinite subset of $\mathcal{F}$, there are $t$ sets that are pierced by a straight line.

\section{Proof of Theorem \ref{exmpl:lines_R3}}
\begin{proof}
  Let $t\in\mathbb{N}$.
  We construct $\mathcal{F}$ in three steps.\\

  \noindent\textit{Step-I:} 
  In this step, we create a sequence of compact sets $\mathcal{K}:=\{K_n\}_{n\in\mathbb{N}}$ in $[0,1]$. \\

  Choose any $\delta\in(0,1)$.
  For every $i\in\mathbb{N}$, let $I^i_1,\dots,I^i_{r_i}$ be an open cover of $[0,1]$ 
  consisting of open intervals of length $\frac{1-\delta}{2^i}$.
  Let $\mathcal{K}_i$ be the family of all sets of the form $[0,1]\setminus (\cup_{j\in[2^i]}I^i_{n_j})$,
  where $n_j\in[r_i]$.
  Then $\mathcal{K}_i$ is a finite collection of compact sets in $[0,1]$, and every set 
  in $\mathcal{K}_i$ has measure at least $\delta$.
  Define $\mathcal{K}=\cup_{i\in\mathbb{N}}\mathcal{K}_i$.\\

  \noindent\textit{Step-II:}
  In this step, we shall show that every infinite subset of $\mathcal{K}$ contains 
  $t$ sets whose intersection is nonempty.\\
  
  Let $\{K'_n\}_{n\in\mathbb{N}}$ be an infinite subset of $\mathcal{K}$.
  Define $S_t := \cup_{\alpha\in\mathbb{N}^t, \alpha_1<\dots<\alpha_t}(\cap_{j\in[t]}K'_{\alpha_j})$, and let 
  $S_i := \cup_{\alpha\in\mathbb{N}^i,\alpha_1<\dots<\alpha_i}(\cap_{j\in[i]}K'_{\alpha_j})\setminus (\cup^t_{s=i+1}S_s)$
  for $i\in\{1,\dots,t-1\}$.
  In other words, $S_t$ constitutes of all points in $\cup_{n\in\mathbb{N}}K'_n$ that lie in 
  $K'_n$ for {\em at least} $t$ distinct values of $n$, and $S_i$ constitutes of all 
  points in $\cup_{n\in\mathbb{N}}K'_n$ that lie in $K'_n$ for {\em exactly} $i$ distinct 
  values of $n\in\mathbb{N}$ when $i\in\{1,\dots,t-1\}$.
  We claim that $\mu(S_t)>0$, where $\mu$ denotes the Lebesgue measure on $\mathbb{R}$.
  If not, then let $\mu(S_t)=0$.
  For each $i\in[t-1]$, every $S_i$ can be written as $\cup_{\alpha\in\mathbb{N}^i,\alpha_1<\dots<\alpha_i} J^i_\alpha$,
  where $J^i_{\alpha}=\cap_{j\in[i]}K'_{\alpha_j}\setminus (\cup_{s=i+1}^t S_s)$.
  In other words, given $\alpha\in\mathbb{N}^i$, $\alpha_1<\dots<\alpha_i$, $J^i_{\alpha}$ contains all points in the intersection of $K'_{\alpha_1},\dots,K'_{\alpha_i}$, which are not in any other $K'_n$, $n\notin\{\alpha_1,\dots,\alpha_i\}$, and $S_i$ is the union of all such $J^i_{\alpha}$.
  Then, $J^i_{\alpha}\cap J^i_{\beta}=\emptyset$ whenever $\alpha\neq\beta$, and also, as $S_i\cap S_{i'}=\emptyset$
  whenever $i\neq i'$, $i'\in[t]$, we have $J^i_{\alpha}\cap J^{i'}_{\beta}=\emptyset$
  whenever $i\neq i'$, $i'\in[t-1]$.

  Therefore, $\mu(S_i)=\sum_{\alpha\in\mathbb{N}^i,\alpha_1<\dots<\alpha_i}\mu(J^i_{\alpha})$.
  Also, $\mu(\cup_{n\in\mathbb{N}}K'_n)=\sum_{i\in[t-1]}\mu(S_i)=s$.
  Then, there is an $m\in\mathbb{N}$ such that $\sum_{i\in[t-1]}(\sum_{\alpha\in[m]^i,\alpha_1<\dots<\alpha_i}\mu(J^i_{\alpha}))>s-\frac{\delta}{2}$.
  But this is not possible, as $K'_{m+1}\subset (\cup_{n\in\mathbb{N}}K'_n)\setminus (\sum_{i\in[t-1]}(\sum_{\alpha\in[m]^i,\alpha_1<\dots<\alpha_i}J^i_{\alpha}))$.
  Therefore, $\mu(S_t)>0$, showing that $S_t$ is non-empty, which implies that $\{K'_n\}_{n\in\mathbb{N}}$ contains $t$ sets that are pierced by a line.\\

  \noindent\textit{Step-III:}
  In this step, we construct $\mathcal{F}$.\\
  
  Denote the rational numbers in $[0,1]$ by $Q_0=\{q_1,q_2,\dots\}$, 
  For every $n\in\mathbb{N}$, let $\mathcal{M}_n$ denote all sets in $\mathcal{K}$ 
  that contains $q_n$ but does not contain $q_1,\dots,q_{n-1}$.
  Clearly, for every $n\in\mathbb{N}$, $\mathcal{M}_n$ is an infinite collection 
  of sets in $\mathcal{K}$, and $\mathcal{K}=\cup_{n\in\mathbb{N}}\mathcal{M}_n$.
  Pick an infinite sequence $Q_1:=\{q^1_n\}_{n\in\mathbb{N}}$ of distinct rational numbers in $[0,1]$ 
  such that $q^1_n\to q_1$ as $n\to\infty$, and $q^1_n\neq q_1$ $\forall n\in\mathbb{N}$.
  Define a bijection $g_1:Q_1\to\mathcal{M}_1$.
  We inductively define $Q_m$, $m\in\{2,3,\dots\}$.
  Pick an infinite sequence $Q_m:=\{q^m_n\}_{n\in\mathbb{N}}$ of distinct 
  rational numbers in $[0,1]\setminus (\cup_{i=1}^{m-1}Q_i)$ 
  such that $q^m_n\to q_m$ as $n\to\infty$. 
  Define a bijection $g_m:Q_m\to\mathcal{M}_m$.

  Let $T=\cup_{m\in\mathbb{N}}Q_m$.
  Note that $T$ is dense in $[0,1]$ for every $m\in\mathbb{N}$.
  Define a bijection $g:T\to\mathcal{K}$ by $g(q)=g_m(q)$ if $q\in Q_m$, $m\in\mathbb{N}$.

  Let, for every $r\in[0,1]$, $l_r$ denote the straight line 
  $x = r$, $z = ry$ in $\mathbb{R}^3$.
  Choose a sequence $\{\epsilon_n\}_{n\in\mathbb{N}}$ of very small positive 
  numbers such that $\epsilon_{n+1}\ll\epsilon_n\ll 1$ for every $n\in\mathbb{N}$, 
  and $\epsilon_n\to 0$ as $n\to\infty$.
  Define a bijection $f:T\to\mathbb{N}$.
  For every $q\in T$, define the plane $\rho_q$ in $\mathbb{R}^3$ by 
  the equation $y=q+\epsilon_{f(q)}x$, and define the 
  set $C_q:=\mathrm{conv}(\rho_q\cap(\cup_{r\in g(q)}l_r))$ (see Figure \ref{fig:Cq}).
  It is easy to see that $C_q$ is compact, and $C_q\subset[0,2]^3$ for every $q\in T$.
  We define $\mathcal{F}$ to be the collection $\{C_q\}_{q\in T}$.
  Clearly, for any $t\in\mathbb{N}$, in every infinite subset of $\mathcal{F}$,
  there are $t$ sets that are pierced by some $l_r$, $r\in[0,1]$.\\

  \begin{figure}[h]
      \centering
      \includegraphics[width=1\linewidth]{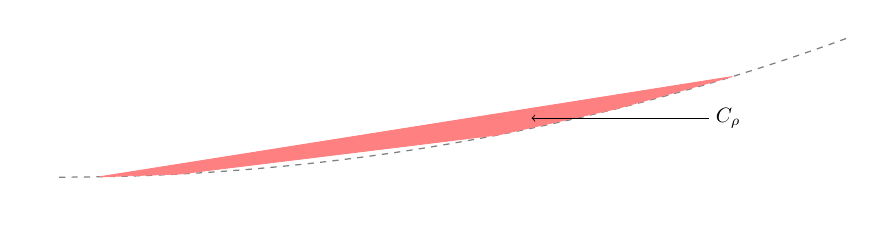}
      \caption{An example of a $C_q$}
      \label{fig:Cq}
  \end{figure}
  
  Now, we show that $\mathcal{F}$ does not have a finite-sized line transversal.
  Note that for every $r\in[0,1]$, $l_r$ lies on the hyperbolic paraboloid $\Sigma$
  given by the equation $z=xy$.
  For every $q\in T$, $\rho_q$ is a plane that is almost perpendicular 
  to the $y$-axis, tilted slightly due to the $\epsilon_{f(q)}$ factor.
  The intersection of $l_r$ and a plane $y=b$ is the point $(r,b,rb)$ that lies on the straight 
  line $y=b$, $z=bx$, which is a line on $\Sigma$.
  So, the intersection that the slightly tilted $\rho_q$ makes with $l_r$ is "very close"
  to the line $y=q$, $z=qx$, due to $\epsilon_{f(q)}$ being very small, 
  and each $C_q$ is a very narrow convex set.
  The intersection of every $\rho_q$ with $\Sigma$ is a parabola in $\rho_q$
  given by the equations $z=q x + \epsilon_{f(q)} x^2$, $y = q + \epsilon_{f(q)} x$,
  and $C_q$ is the convex hull of a compact set on the parabola.
  We can show that the vertical distance of any point on $C_q$ from $\Sigma$
  is bounded above by $\epsilon_{f(q)}$.

  Every straight line in $\mathbb{R}^3$ that intersects $\Sigma$ either lies on $\Sigma$, or intersects 
  $\Sigma$ in at most two points, and the lines that lie on $\Sigma$ are of the form 
  $z=bx$, $y=b$, or $z=cy$, $x=c$.
  Let, if possible, $\mathcal{L}$ be a finite collection of lines that 
  is a transversal for $\mathcal{F}$. 
  We divide $\mathcal{L}$ into three collections $\mathcal{L}_1,\mathcal{L}_2$ and $\mathcal{L}_3$,
  where $\mathcal{L}_1$ contains all lines of the form $z=cy$, $x=c$,
  $\mathcal{L}_2$ contains all lines of the form $z=bx$, $y=b$,
  and $\mathcal{L}_3=\mathcal{L}\setminus (\mathcal{L}_1\cup\mathcal{L}_2)$.

  Every line in $\mathcal{L}_1$ must be an $l_r$ for some $r\in[0,1]$.
  Let $\mathcal{L}_1=\{l_{r_1},\dots,l_{r_m}\}$.
  Note that if $r\notin g(q)$ for some $q\in T$, then $l_r$ does not 
  intersect $C_q$.
  That is because if $a$ is a point on the convex parabolic arc $A_q\subset \rho_q\cap\Sigma$ between the points
  $(0,q,0)$ and $(1,q+\epsilon_{f(q)},q+\epsilon_{f(q)})$ 
  on which $\rho_q\cap (\cup_{r\in g(q)}l_r)$ lies, 
  then $a$ does not lie in the convex hull of any subset of $A_q$ that 
  does not contain $a$.
  So, a line in $\mathcal{L}_1$ intersects $C_q$ if and only if $r_i\in g(q)$ for some $i\in[m]$.
  Let $T_1$ be the subset of $T$ such that $q\in T_1$ if and only if $q\in T$ and $g(q)$
  does not contain $r_i$ for every $i\in[m]$.
  We claim that there is an $n'\in\mathbb{N}$ such that for any $n\in\mathbb{N}$, $n\geq n'$, $T_1\cap Q_n$ contains infinitely many points.
  There is an $n'\in\mathbb{N}$ such that $\{q_{n'},q_{n'+1},\dots\}\cap\{r_1,\dots,r_m\}=\emptyset$.
  Let, without loss of generality, $r_1,\dots,r_{m'}$, $m'\leq m$,
  be irrational.
  We show that for every $n\in\mathbb{N}$, there are infinitely many 
  sets in $g(Q_n)$ that do not contain $r_1,\dots,r_{m'}$, which will prove the claim. 
  Pick any $n\in\mathbb{N}$ and let $N\in\mathbb{N}$ be a number such that $2^N>m'+n$, $\frac{1-\delta}{2^N}\ll 
  \mathrm{dist}(r_i,q_n)$ $\forall i\in[m']$,
  and $\frac{1-\delta}{2^N}\ll \mathrm{dist}(q_k,q_n)$ $\forall k\in[n-1]$.
  Therefore, in any finite cover of $[0,1]$ by open intervals of length $\frac{1-\delta}{2^i}$
  for $i\geq N$, $q_n$ does not lie in any interval containing any of 
  the points $r_1,\dots,r_{m'},q_1,\dots,q_{n-1}$.
  As the sets in $\mathcal{K}_i$, $i\geq N$, are obtained by removing every collection of $2^i$ intervals from $[0,1]$ picked from an open cover of $[0,1]$, and $2^i>n+m'$, there are infinitely many sets in $\mathcal{M}_n$ that do not contain 
  $r_1,\dots,r_{m'}$.

  Let $\mathcal{L}_2=\{l^2_1,\dots,l^2_{m_2}\}$, where $l^2_j$ is given by the equations 
  $y=r_j$, $z=r_j x$, $r_j\in [0,1]$ for all $j\in[m_2]$.
  Let $T_2$ be the subset of $T_1$ such that $q\in T_2$ if and only if 
  $q\in T_1$ and $[0,2]^3\cap \rho_q$ does not intersect the plane $y=r_j$
  for any $j\in[m_2]$.
  Note that since for every $q\in T$, $C_q\subset[0,2]^3\cap\rho_q$, if $[0,2]^3\cap\rho_q$
  does not intersect the plane $y=r_j$, then $l^2_j$ does not intersect $C_q$.
  We claim that there is a $n''\in\mathbb{N}$ such that for every $n\in\mathbb{N}$, $n\geq n''$, $Q_n\cap T_2$ is an infinite set.
  There is an $n''\in\mathbb{N}$ with $n''\geq n'$ such that $\{q_{n''},q_{n''+1},\dots\}\cap \{r_1,\dots,r_{m_2}\}=\emptyset$.
  Let $n\in\mathbb{N}$, $n\geq n''$.
  Then, there is an $\epsilon'_n>0$ such that $|r_j - q_n|>2\epsilon'_n$ for all $j\in[m_2]$.
  As $T_1$ contains infinitely many points from $Q_n$, there are infinitely many points $q$ in $T_1\cap Q_n$ such that $|q-q_n|<\epsilon'_n$, and the value of $\epsilon_{f(q)}$ is so small that the plane $\rho_q$ has a very small tilt, and $\rho_q\cap [0,2]^3$ does not intersect the plane $y=r_j$ for any $j\in[m_2]$.
  
  Let $\mathcal{L}_3:=\{l^3_1,\dots,l^3_{m_3}\}$.
  If $l^3_j\in\mathcal{L}_3$, then $l^3_j$ intersects $\Sigma$ 
  in at most two points, say $p^{j}_1$ and $p^{j}_2$.
  Choose positive numbers $\xi_j$ such that $4\sum_{j\in[m_3]}\xi_j\ll 1$.
  For each $j\in[m_3]$, let $B^j_1$ and $B^j_2$ denote two open balls with radius 
  length $\xi_j$, and and centers at $p^j_1$ and $p^j_2$, respectively.
  There is an $\epsilon >0$ for which the vertical distance between any point on $l^3_j$ and 
  $(\Sigma\cap [0,2]^3)\setminus \cup_{k\in[m_3]}(B^k_1\cup B^k_2)$
  is at least $\epsilon$ for every $j\in[m_3]$ (if a vertical line through a point 
  on $l^3_j$ does not intersect $(\Sigma\cap [0,2]^3)\setminus \cup_{k\in[m_3]}(B^k_1\cup B^k_2)$,
  we take the distance to be $\infty$).
  Since $4\sum_{j\in[m_3]}\xi_j\ll 1$, 
  there are infinitely many $q\in T_2$ for which 
  $\rho_q$ does not intersect $B^j_i$ for any $i\in\{1,2\}$ and any $j\in[m_3]$.
  Since the vertical distance between $C_q$ and $\Sigma$ is bounded above 
  by $\epsilon_{f(q)}$ for every $q\in T_2$, there are infinitely many $q\in T_2$
  for which $C_q$ does not intersect $l^3_j$ for any $j\in[m_3]$.

  Therefore, $\mathcal{L}$ cannot be a line transversal for $\mathcal{F}$, which 
  is a contradiction.
  This concludes the proof.
\end{proof}

\bibliographystyle{alpha}
\bibliography{references}

\begin{thebibliography}{AKMM02}

\bibitem[AK92]{alon1992piercing}
Noga Alon and Daniel~J Kleitman.
\newblock Piercing convex sets and the hadwiger-debrunner (p, q)-problem.
\newblock {\em Advances in Mathematics}, 96(1):103--112, 1992.

\bibitem[AK95]{alon1995bounding}
Noga Alon and Gil Kalai.
\newblock Bounding the piercing number.
\newblock {\em Discrete \& Computational Geometry}, 13(3):245--256, 1995.

\bibitem[AKMM02]{alon2002transversal}
Noga Alon, Gil Kalai, Jiri Matousek, and Roy Meshulam.
\newblock Transversal numbers for hypergraphs arising in geometry.
\newblock {\em Advances in Applied Mathematics}, 29(1):79--101, 2002.

\bibitem[CGH24]{cheong2024some}
Otfried Cheong, Xavier Goaoc, and Andreas~F Holmsen.
\newblock Some new results on geometric transversals.
\newblock {\em Discrete \& Computational Geometry}, 72(2):674--703, 2024.

\bibitem[CGN25a]{chakraborty2025finitektransversalsinfinitefamilies}
Sutanoya Chakraborty, Arijit Ghosh, and Soumi Nandi.
\newblock Finite k-transversals of infinite families of fat convex sets, 2025.

\bibitem[CGN25b]{chakraborty2025stabbing}
Sutanoya Chakraborty, Arijit Ghosh, and Soumi Nandi.
\newblock Stabbing boxes with finitely many axis-parallel lines and flats.
\newblock {\em Discrete Mathematics}, 348(2):114269, 2025.

\bibitem[HM04]{holmsen2004no}
Andreas Holmsen and Jiri Matousek.
\newblock No helly theorem for stabbing translates by lines in r3.
\newblock {\em Discrete \& Computational Geometry}, 31(3):405--410, 2004.

\bibitem[JP24]{jung2024note}
Attila Jung and D{\"o}m{\"o}t{\"o}r P{\'a}lv{\"o}lgyi.
\newblock A note on infinite versions of $(p, q) $-theorems.
\newblock {\em arXiv preprint arXiv:2412.04066}, 2024.

\bibitem[KP22]{keller2022}
Chaya Keller and Micha~A Perles.
\newblock An ($\aleph_0$, k+ 2)-theorem for k-transversals.
\newblock In {\em 38th International Symposium on Computational Geometry (SoCG
  2022)}, pages 50--1. Schloss Dagstuhl--Leibniz-Zentrum f{\"u}r Informatik,
  2022.

\bibitem[KP23]{keller2023aleph_0}
Chaya Keller and Micha~A Perles.
\newblock An $(\aleph_0, k+ 2) $-theorem for $ k $-transversals.
\newblock {\em arXiv preprint arXiv:2306.02181}, 2023.

\end{thebibliography}

\end{document}